\documentclass[12pt, reqno]{article}

\usepackage[colorlinks=true,
linkcolor=webgreen,
filecolor=webbrown,
citecolor=webgreen]{hyperref}
\usepackage{amsthm}
\usepackage{graphics, amsmath, amssymb}
\usepackage{fullpage}
\usepackage{mathtools} 
\usepackage{tikz}
\usepackage[numbers]{natbib}
\usepackage{listings}
\usepackage{enumitem}

\definecolor{webgreen}{rgb}{0,.5,0}
\definecolor{webbrown}{rgb}{.6,0,0}

\setlist{  
  listparindent=\parindent,
  parsep=0pt,
}
\usetikzlibrary{arrows}

\newcommand{\EQ}[1]{
\begin{equation}\begin{split}#1\end{split}\end{equation}
}
\newcommand{\comment}[1]{}
\newcommand{\Ker}{{\operatorname{supp}}}
\newcommand{\pop}{{s_2}}
\newcommand{\Char}{{\operatorname{char}}}
\newcommand{\seqnum}[1]{\href{http://oeis.org/#1}{\underline{#1}}}

\title{Surjectivity of a Binary Analog of the Enots Wolley Sequence}

\theoremstyle{plain}
\newtheorem{lemma}{Proposition}
\newtheorem{lemma1}[lemma]{Lemma}
\newtheorem{lemma2}[lemma]{Lemma}

\theoremstyle{plain}
\newtheorem{prop}{Proposition}
\newtheorem{theorem1}[prop]{Theorem}
\newtheorem{theorem2}[prop]{Theorem}
\newtheorem{theorem3}[prop]{Theorem}

\newtheorem{theorem6}[prop]{Theorem}
\newtheorem{theorem7}[prop]{Theorem}

\theoremstyle{definition}
\newtheorem{defn}{Definition}
\newtheorem{definition1}[defn]{Definition}
\newtheorem{definition2}[defn]{Definition}
\newtheorem{definition3}[defn]{Definition}
\newtheorem{definition4}[defn]{Definition}
\newtheorem{definition5}[defn]{Definition}
\newtheorem{definition6}[defn]{Definition}
\newtheorem{definition7}[defn]{Definition}

\begin{document}

\begin{center}
\vskip 1cm{\LARGE\bf 
The Binary Enots Wolley Sequence\\
}
\vskip 1cm
\large
Nathan Nichols\\
Milwaukee, Wisconsin \\ 
United States \\
\href{mailto:nathannichols454@gmail.com}{\tt nathannichols454@gmail.com} \\
\end{center}

\begin{abstract}
It is an open conjecture that the Enots Wolley sequence is surjective onto the set of positive integers with a binary weight of at least 2. In this paper, this property is proved for an analog of the Enots Wolley sequence which operates on the binary representation of a number rather than the prime factorization.
\end{abstract}

\section{Introduction}
The Yellowstone permutation is sequence \seqnum{A098550} in the OEIS. It was defined and studied by Applegate et al.\ who proved the key property that every positive integer appears as a term \cite{YSPerm, YellowstoneOEIS}. The Enots Wolley sequence \seqnum{A336957} is closely related to the Yellowstone permutation. It is still an open question whether or not every eligible number appears in the Enots Wolley sequence \cite{EnotsOEIS, SloaneTalk}. 

This paper is concerned with a variation of the Enots Wolley sequence which operates on the binary representation of a number rather than the prime factorization. This variation (sequence \seqnum{A338833}) is called the \textit{Binary} or \textit{Set Theoretic} Enots Wolley sequence \cite{BinaryEnotsOEIS}.  For this variation, the property that every number with a binary weight of at least two appears in the sequence is proved. The exact definition of the Binary Enots Wolley sequence is as follows:

\begin{definition1}
For $n\geq 1$, the Binary Enots Wolley sequence $a(n)$ is defined as the lexicographically earliest infinite sequence of distinct positive integers with the property that $a(n)$ is the any number $m$ not yet in the sequence such that the binary expansions of $m$ and $a(n-1)$ have a $1$ in the same position, but the positions of the $1$s in the binary expansions of $m$ and $a(n-2)$ are disjoint.
\end{definition1}

The remaining definitions and theorems in this section have direct analogs for the Enots Wolley sequence \cite{EnotsOEIS}.

\begin{definition2} Let $k$ be a positive integer. Let $\Ker(k)$, the \textit{support} of $k$, denote the set of positions of $1$s in the binary expansion of $k$. Thus $\Ker(15) = \{0,1,2,3\}$ and $\Ker(1) = \{0\}$.
\end{definition2}

\begin{definition3}
For a positive integer $k$, the \textit{characteristic function} $\Char_k(i): \mathbb{N}\to \{0, 1\}$ is 
\EQ{
\Char_k(i) := \begin{cases} 
      1 & \Ker(k) \cap \Ker(a(i)) \neq \emptyset, \\
      0 & \textrm{otherwise}.\\
\end{cases}
}
 
\end{definition3}

\begin{definition4}
The \textit{truth table} of $\Char_k(i)$ is the sequence $(\Char_k(i))_{i>0}$.
\end{definition4}

\begin{definition5} For a positive integer $k$, let $\pop(k)$ denote the number of $1$s in the binary expansion of $k$ (the \textit{binary weight} of $k$.)
\end{definition5}

The following theorem serves as a more useful form of Definition 1:

\begin{theorem1} For $n > 2$, $a(n)$ is the smallest number $m$ not yet in the sequence such that: 
\begin{enumerate}
\item[i)] $\Ker(m)\cap \Ker(a(n-1))$ is nonempty,
\item[ii)] $\Ker(m) \cap \Ker(a(n-2))$ is empty, and
\item[iii)] the set $\Ker(m)\setminus \Ker(a(n-1))$ is nonempty.
\end{enumerate}
\end{theorem1}
\begin{proof}
Conditions (i) and (ii) are already part of Definition 1. Condition (iii) is a consequence of (i) and (ii) because if (iii) fails for some term $a(n)$, there is no way for property (ii) to hold for the term $a(n+1)$. This proves that if the Binary Enots Wolley exists, then it must satisfy properties (i)-(iii). 

The Binary Enots Wolley sequence is defined in Definition 1 as the lexicographically earliest sequence in some set  $S$ of sequences. To prove the Binary Enots Wolley sequence exists, it suffices to show that the set $S$ is not empty. Here is one example of a sequence of binary numbers $e_i\in S$:
\begingroup
\setlength\arraycolsep{0pt}
\[
\begin{matrix*}[l]
	&e_1=1111 \\
	&e_2=001111 \\
	&e_3=00001111\\
	&e_4=0000001111\\
	&\vdots
	
\end{matrix*}
\]
\endgroup 

It remains to be shown that $a(n)$ is the \textit{smallest} number $m$ satisfying properties (i)-(iii). If it were possible to construct a sequence in $S$ term-by-term by always choosing $a(n)$ to be the least number $m$ satisfying conditions (i)-(iii), then this ``greedy" strategy would produce the lexicographically earliest sequence in $S$. Hence, it suffices to show that this greedy strategy does not lead to any dead ends. For this, observe that there will always be at least one possible choice for $a(n)$ if $a(n-1)$ and $a(n-2)$ satisfy conditions (i)-(iii).
\end{proof}

An immediate implication of Theorem 1 is that only numbers with a binary weight of at least $2$ can appear in the sequence. See Figure 1 for a table of the first terms of the Enots Wolley and Binary Enots Wolley sequences. 

\begin{definition6}
Let $m$ be a positive integer such that $\pop(m)\geq 2$. Say that $m$ is a \textit{candidate} for the $n$th term if properties (i), (ii), and (iii) of Theorem 1 hold.
\end{definition6}

\newpage
\begin{theorem7}
If $n$ is the least positive integer such that $v\in \Ker(a(n))$, the term $a(n)$ has a binary weight of $2$.
\end{theorem7}\begin{proof}
Suppose $n$ is the least positive integer such that $v \in \Ker(a(n))$ and let $S$ be the set of all candidates for the $n$th term. By Theorem 1, the term $a(n)$ is the least element of $S$ which does not appear in the first $n-1$ terms of the sequence. 

Define $T$ to be the set of elements of $\Ker(a(n-1))$ that are not also contained in $\Ker(a(n-2))$. By Theorem 1, for every $m\in S$ there is at least one element of $T$ contained in $\Ker(m)$. The least element of $S$ involving $2^v$ is $2^v+2^w$ where $w$ is the least element of $T$. Because $2^v$ does not appear in the first $n-1$ terms of the sequence, neither does the number $2^v+2^w$. Thus, $a(n)=2^v+2^w$.
\end{proof}

\begin{definition7}
If $a(n)=2^v + 2^w$ is the first term of the sequence such that $v\in \Ker(a(n))$, say that $2^v$ is \textit{introduced} by $2^w$.  
\end{definition7}

\begin{figure}
\centering
\input{enots.tikz}
\begin{tikzpicture}[scale=0.54]
   \begin{scope}
      \draw[step=1cm] (0, 0) grid (8, 35);
\fill[orange!20] (0,34) rectangle (1,35);
\fill[orange!20] (1,34) rectangle (2,35);
      \node[anchor=center,font=\small] at (1.5, 34.5){1};
\fill[orange!20] (2,34) rectangle (3,35);
      \node[anchor=center,font=\small] at (2.5, 34.5){2};
\fill[orange!20] (3,34) rectangle (4,35);
      \node[anchor=center,font=\small] at (3.5, 34.5){4};
\fill[orange!20] (4,34) rectangle (5,35);
      \node[anchor=center,font=\small] at (4.5, 34.5){8};
\fill[orange!20] (5,34) rectangle (6,35);
      \node[anchor=center,font=\small] at (5.5, 34.5){16};
\fill[orange!20] (6,34) rectangle (7,35);
      \node[anchor=center,font=\small] at (6.5, 34.5){32};
\fill[orange!20] (7,34) rectangle (8,35);
      \node[anchor=center,font=\small] at (7.5, 34.5){64};
\fill[orange!20] (0,33) rectangle (1,34);
      \node[anchor=center,font=\small] at (0.5, 33.5){1};
      \node[anchor=center] at (1.5, 33.5){1};
\fill[orange!20] (0,32) rectangle (1,33);
      \node[anchor=center,font=\small] at (0.5, 32.5){2};
      \node[anchor=center] at (1.5, 32.5){1};
      \node[anchor=center] at (2.5, 32.5){1};
\fill[orange!20] (0,31) rectangle (1,32);
      \node[anchor=center,font=\small] at (0.5, 31.5){3};
      \node[anchor=center] at (2.5, 31.5){1};
      \node[anchor=center] at (3.5, 31.5){1};
\fill[orange!20] (0,30) rectangle (1,31);
      \node[anchor=center,font=\small] at (0.5, 30.5){4};
      \node[anchor=center] at (3.5, 30.5){1};
      \node[anchor=center] at (4.5, 30.5){1};
\fill[orange!20] (0,29) rectangle (1,30);
      \node[anchor=center,font=\small] at (0.5, 29.5){5};
      \node[anchor=center] at (1.5, 29.5){1};
      \node[anchor=center] at (4.5, 29.5){1};
\fill[orange!20] (0,28) rectangle (1,29);
      \node[anchor=center,font=\small] at (0.5, 28.5){6};
      \node[anchor=center] at (1.5, 28.5){1};
      \node[anchor=center] at (5.5, 28.5){1};
\fill[orange!20] (0,27) rectangle (1,28);
      \node[anchor=center,font=\small] at (0.5, 27.5){7};
      \node[anchor=center] at (2.5, 27.5){1};
      \node[anchor=center] at (5.5, 27.5){1};
\fill[orange!20] (0,26) rectangle (1,27);
      \node[anchor=center,font=\small] at (0.5, 26.5){8};
      \node[anchor=center] at (2.5, 26.5){1};
      \node[anchor=center] at (4.5, 26.5){1};
\fill[orange!20] (0,25) rectangle (1,26);
      \node[anchor=center,font=\small] at (0.5, 25.5){9};
      \node[anchor=center] at (1.5, 25.5){1};
      \node[anchor=center] at (3.5, 25.5){1};
      \node[anchor=center] at (4.5, 25.5){1};
\fill[orange!20] (0,24) rectangle (1,25);
      \node[anchor=center,font=\small] at (0.5, 24.5){10};
      \node[anchor=center] at (3.5, 24.5){1};
      \node[anchor=center] at (5.5, 24.5){1};
\fill[orange!20] (0,23) rectangle (1,24);
      \node[anchor=center,font=\small] at (0.5, 23.5){11};
      \node[anchor=center] at (5.5, 23.5){1};
      \node[anchor=center] at (6.5, 23.5){1};
\fill[orange!20] (0,22) rectangle (1,23);
      \node[anchor=center,font=\small] at (0.5, 22.5){12};
      \node[anchor=center] at (1.5, 22.5){1};
      \node[anchor=center] at (6.5, 22.5){1};
\fill[orange!20] (0,21) rectangle (1,22);
      \node[anchor=center,font=\small] at (0.5, 21.5){13};
      \node[anchor=center] at (1.5, 21.5){1};
      \node[anchor=center] at (3.5, 21.5){1};
\fill[orange!20] (0,20) rectangle (1,21);
      \node[anchor=center,font=\small] at (0.5, 20.5){14};
      \node[anchor=center] at (2.5, 20.5){1};
      \node[anchor=center] at (3.5, 20.5){1};
      \node[anchor=center] at (4.5, 20.5){1};
\fill[orange!20] (0,19) rectangle (1,20);
      \node[anchor=center,font=\small] at (0.5, 19.5){15};
      \node[anchor=center] at (4.5, 19.5){1};
      \node[anchor=center] at (5.5, 19.5){1};
\fill[orange!20] (0,18) rectangle (1,19);
      \node[anchor=center,font=\small] at (0.5, 18.5){16};
      \node[anchor=center] at (1.5, 18.5){1};
      \node[anchor=center] at (5.5, 18.5){1};
      \node[anchor=center] at (6.5, 18.5){1};
\fill[orange!20] (0,17) rectangle (1,18);
      \node[anchor=center,font=\small] at (0.5, 17.5){17};
      \node[anchor=center] at (1.5, 17.5){1};
      \node[anchor=center] at (2.5, 17.5){1};
      \node[anchor=center] at (3.5, 17.5){1};
\fill[orange!20] (0,16) rectangle (1,17);
      \node[anchor=center,font=\small] at (0.5, 16.5){18};
      \node[anchor=center] at (2.5, 16.5){1};
      \node[anchor=center] at (7.5, 16.5){1};
\fill[orange!20] (0,15) rectangle (1,16);
      \node[anchor=center,font=\small] at (0.5, 15.5){19};
      \node[anchor=center] at (4.5, 15.5){1};
      \node[anchor=center] at (7.5, 15.5){1};
\fill[orange!20] (0,14) rectangle (1,15);
      \node[anchor=center,font=\small] at (0.5, 14.5){20};
      \node[anchor=center] at (1.5, 14.5){1};
      \node[anchor=center] at (4.5, 14.5){1};
      \node[anchor=center] at (5.5, 14.5){1};
\fill[orange!20] (0,13) rectangle (1,14);
      \node[anchor=center,font=\small] at (0.5, 13.5){21};
      \node[anchor=center] at (1.5, 13.5){1};
      \node[anchor=center] at (2.5, 13.5){1};
      \node[anchor=center] at (5.5, 13.5){1};
\fill[orange!20] (0,12) rectangle (1,13);
      \node[anchor=center,font=\small] at (0.5, 12.5){22};
      \node[anchor=center] at (2.5, 12.5){1};
      \node[anchor=center] at (6.5, 12.5){1};
\fill[orange!20] (0,11) rectangle (1,12);
      \node[anchor=center,font=\small] at (0.5, 11.5){23};
      \node[anchor=center] at (3.5, 11.5){1};
      \node[anchor=center] at (6.5, 11.5){1};
\fill[orange!20] (0,10) rectangle (1,11);
      \node[anchor=center,font=\small] at (0.5, 10.5){24};
      \node[anchor=center] at (1.5, 10.5){1};
      \node[anchor=center] at (3.5, 10.5){1};
      \node[anchor=center] at (5.5, 10.5){1};
\fill[orange!20] (0,9) rectangle (1,10);
      \node[anchor=center,font=\small] at (0.5, 9.5){25};
      \node[anchor=center] at (1.5, 9.5){1};
      \node[anchor=center] at (2.5, 9.5){1};
      \node[anchor=center] at (4.5, 9.5){1};
\fill[orange!20] (0,8) rectangle (1,9);
      \node[anchor=center,font=\small] at (0.5, 8.5){26};
      \node[anchor=center] at (4.5, 8.5){1};
      \node[anchor=center] at (6.5, 8.5){1};
\fill[orange!20] (0,7) rectangle (1,8);
      \node[anchor=center,font=\small] at (0.5, 7.5){27};
      \node[anchor=center] at (3.5, 7.5){1};
      \node[anchor=center] at (5.5, 7.5){1};
      \node[anchor=center] at (6.5, 7.5){1};
\fill[orange!20] (0,6) rectangle (1,7);
      \node[anchor=center,font=\small] at (0.5, 6.5){28};
      \node[anchor=center] at (2.5, 6.5){1};
      \node[anchor=center] at (3.5, 6.5){1};
      \node[anchor=center] at (5.5, 6.5){1};
\fill[orange!20] (0,5) rectangle (1,6);
      \node[anchor=center,font=\small] at (0.5, 5.5){29};
      \node[anchor=center] at (1.5, 5.5){1};
      \node[anchor=center] at (2.5, 5.5){1};
      \node[anchor=center] at (7.5, 5.5){1};
\fill[orange!20] (0,4) rectangle (1,5);
      \node[anchor=center,font=\small] at (0.5, 4.5){30};
      \node[anchor=center] at (1.5, 4.5){1};
      \node[anchor=center] at (4.5, 4.5){1};
      \node[anchor=center] at (6.5, 4.5){1};
\fill[orange!20] (0,3) rectangle (1,4);
      \node[anchor=center,font=\small] at (0.5, 3.5){31};
      \node[anchor=center] at (3.5, 3.5){1};
      \node[anchor=center] at (4.5, 3.5){1};
      \node[anchor=center] at (5.5, 3.5){1};
\fill[orange!20] (0,2) rectangle (1,3);
      \node[anchor=center,font=\small] at (0.5, 2.5){32};
      \node[anchor=center] at (3.5, 2.5){1};
      \node[anchor=center] at (7.5, 2.5){1};
\fill[orange!20] (0,1) rectangle (1,2);
      \node[anchor=center,font=\small] at (0.5, 1.5){33};
      \node[anchor=center] at (1.5, 1.5){1};
      \node[anchor=center] at (7.5, 1.5){1};
\fill[orange!20] (0,0) rectangle (1,1);
      \node[anchor=center,font=\small] at (0.5, 0.5){34};
      \node[anchor=center] at (1.5, 0.5){1};
      \node[anchor=center] at (2.5, 0.5){1};
      \node[anchor=center] at (4.5, 0.5){1};
      \node[anchor=center] at (5.5, 0.5){1};
      \node[anchor=center] at (4.0, -0.5) {Binary Enots Wolley};
   \end{scope}
\end{tikzpicture}
\caption{The prime factorizations of the first 34 terms of the Enots Wolley sequence and the binary expansions of the first 34 terms of the binary Enots Wolley sequence. A blank square indicates a multiplicity of zero or a coefficient of zero (respectively.)}
\end{figure}

\newpage
\section{Surjectivity of the Binary Enots Wolley sequence}

\begin{theorem2} Let $k$ be a positive integer such that $\pop(k)\geq 2$. Let $S(k)$ be a set of positive integers $i$ such that $k$ is a candidate for $a(i)$. If there are $k$ or more elements of $S$, then $k$ must appear in the sequence. 
\end{theorem2}
\begin{proof}
If $k$ is a candidate for the $i$th term but $a(i)\neq k$, it must be that $a(i)=k'$ where $k'$ is another candidate for the $i$th term that is less than $k$. Because there are only $k-1$ positive integers $k'$ less than $k$ and each can only appear in the sequence once, this situation can only occur for at most $k-1$ choices of $i\in S(k)$ before $k$ must be the least candidate that has not already occurred in the sequence.
\end{proof}

\begin{lemma2}
For all non-negative integers $m$, the truth table of $\Char_{2^m}$ has infinitely many occurrences of $01$.
\end{lemma2}
\begin{proof}
Since the truth table of a power of two cannot have any occurrences of $111$, it must only be shown that there is no $i'$ such that $i > i'$ implies $\Char_{2^m}(i)=0$. Assume that this is the case. Let $2^a$ be the highest bit that has occurred in the sequence before the $(i')$th term. Let $2^b$ be a sufficiently large power of two so that $b-a>N$ where $N$ is some large number (for instance, $N=100$.) Let $n$ be the least positive integer such that $b\in \Ker(a(n))$. It is now necessary to make the following three definitions:
\begin{enumerate}
\item[-] $U :=$ The bit that introduces $2^b$ in the $n$th term.
\item[-] $V :=$ A power of $2$ distinct from $U$ such that $2^a< V < 2^b$.
\item[-]   $W := 2^m$.
\end{enumerate}
It will be shown that unless all of the $N-1$ choices of $V$ are present in $a(n-2)$, the number $U+V+W$ is a smaller unused candidate for $a(n)$ than $U+2^b$. For one, the number $U+V+W$ is less than $2^b + U$ because $V + W < 2^b$. The other parts of this claim follow from the following properties of the three bits $U$, $V$ and $W$:
\begin{enumerate}
\item[- $U$:]  Ensures properties (i) and (iii) for $U+V+W$ to be a candidate for the $n$th term hold.
\item[- $V$:] Because $W$ has never occurred together with $V$, this ensures $U+V+W$ has not appeared in the first $n-1$ terms of the sequence.
\item[- $W$:]  Ensures that if $a(n)=U+V+W$, then $\Char_{2^m}(n) = 1$ (contradicting the assumption that $\Char_{2^m}(n) = 0$.)
\end{enumerate}

This proves that as long as $V$ isn't barred from appearing in the term $a(n)$ by property (ii), the number $U+V+W$ is a smaller candidate for $a(n)$ that does not appear in the first $n-1$ terms of the sequence. In the case that all possible choices of $V$ are present in $a(n-2)$, it is possible to construct an alternative value of $a(n-2)$ that produces a lexicographically earlier sequence as follows: 

Because $b-a<N$, there are $N-1$ possible choices for $V$. Hence, if all possible choices of $V$ are present in $a(n-2)$, then $a(n-2)$ must have a binary weight of at least $N-1$. Since none of the bits that are possible choices for $V$ have yet appeared in the sequence in a term together with $2^m$, it is possible to erase one of the bits that are possible choices for $V$ and replace it with $2^m$ in $a(n-2)$ to obtain a smaller unused candidate for the $(n-2)$nd term. 
\end{proof}

\begin{lemma1}
Let $k=2^p + 2^q$ be a positive integer with $\pop(k)=2$. If there are an infinite number of occurrences of $01$ in the truth table of $\Char_k(i)$, then $k$ must appear in the sequence.
\end{lemma1}
\begin{proof}
Let $S$ be an infinite set of positive integers such that $\Char_k(j)=0$ and $\Char_k(j+1)=1$ for all $j\in S$. For each $j \in S$, properties (i) and (ii) of Theorem 1 always hold for $k$ to be a candidate for the $(j+2)$nd term. If property (iii) for $k$ to be a candidate for the $(j+2)$nd term also happens to hold for an infinite number of positions $j\in S$, then $k$ is a candidate an infinite number of times and therefore must appear in the sequence by Theorem 3. In this case, there is nothing to prove.

If property (iii) for $k$ to be a candidate for the $(j+2)$nd term holds only for a finite number of $j\in S$, it must be that $\Ker(k)$ is a subset of $\Ker(a(j+1))$ an infinite number of times. Thus, all but finitely many of the $j \in S$ correspond to instances of the following pattern, where each of the $m_0,m_1,m_2,m_3$ have supports disjoint with $\Ker(k)$ and $d$ is some positive integer such that $\Ker(d)$ is a (possibly empty) subset of $\Ker(k)$:
\EQ{
\begin{matrix}
	&\ 	&(A) &\  &(B) &\ &(C) &\ &(D) &\ \\
\cdots &\to &m_0 &\to &m_1 +k &\to &m_2+d &\to &m_3 &\to\  \cdots \\
\end{matrix}
}
This means that almost every time there is a $1$ after a run of $0$s in the truth table of $\Char_k(i)$, all bits of $k$ must be present in the term corresponding to the leading $1$ (i.e., in position $(B)$.) Note that the term in position $(D)$ must have  support disjoint with $\Ker(k)$. This implies that there can be almost no occurrences of $111$ in truth table of $\Char_k(i)$.

Let $P_i$ and $Q_i$ denote the set of distinct positive integers appearing in the first $i$ terms of the sequence that contain $p$ (but not $q$) and $q$ (but not $p$) in their supports (respectively.) Define $P(i) := |P_i|$ and $Q(i) := |Q_i|$ to be the sizes of these sets. Likewise, let $PQ(i) := |PQ_i|$ denote the number of distinct positive integers appearing in the first $i$ terms of the sequence which contain both $p$ and $q$ in their supports.

The pattern (2) implies that after a finite number of terms, at least half of the new $1$s in the truth table of $\Char_k$ contribute to $PQ_i$, and the rest of the new $1$s contribute to either $P_i$ or $Q_i$. Because every bit must appear infinitely many times by Lemma 1, it must be that there exists $i'$ such that $i > i'$ implies $PQ(i) > P(i)$ or $PQ(i) > Q(i)$. In other words, the function $PQ(i)$ grows faster than at least one of $P(i)$ or $Q(i)$.

On the other hand, every time pattern (2) occurs, both $m_1+2^p$ and $m_1+2^q$ are smaller candidates than $m_1 + k$ for the term in position $(B)$. If the smallest unused candidate for position $(B)$ happens to be $m_1+k$, then it must be the case that both $m_1 + 2^p$ and $m_1 + 2^q$ have already occurred in the sequence. This implies there exists $i'$ such that $i > i'$ implies $Q(i) \geq PQ(i)$ and $P(i) \geq PQ(i)$, which is a contradiction of the other lower bound for $PQ(i)$ from the previous paragraph. Hence, our original assumption that property (iii) for $k$ to be a candidate for the $(j+2)$nd term holds only a finite number of times after an occurrence of $01$ in the truth table of $\Char_k$ is false. By Theorem 3, this implies that $k$ must occur in the sequence. 
\end{proof}

\begin{theorem3} Let $k$ be a positive integer such that $\pop(k) \geq 2$. If there are an infinite number of occurrences of $01$ in the truth table of $\Char_k(i)$, then $k$ must appear in the sequence.
\end{theorem3}
\begin{proof} (Outline.)
In the case that $\pop(k) > 2$, the proof of Lemma 2 breaks down in the fourth paragraph (starting with ``Let $P_i$ and $Q_i$...") The proof from the fourth paragraph onward can be adapted to the case that $\pop(k)> 2$ by observing that every $m_1+s$ where $\Ker(s)$ is a non-empty proper subset of $\Ker(k)$ is a smaller candidate for position $(B)$. Thus, occurrences of terms whose support contains $\Ker(k)$ must be even more rare than in the case of $\pop(k)=2$. Yet, the proof of Lemma 2 still leads to the situation where at least half of all but finitely many $1$s in the truth table of $\Char_k$ correspond to terms that contain all of $\Ker(k)$ in their support.
\end{proof}

\begin{lemma2} Let $k =2^p + 2^q$ be a positive integer with $\pop(k) = 2$. Then, there is no positive integer $i'$ such that $i>i'$ implies $\Char_k(i)=1$.
\end{lemma2}
\begin{proof}

Suppose that there exists $i'$ such that $i > i'$ implies $\Char_{k}(i)=1$. Using Iverson bracket notation, define a function $f:\mathbb{N}_{>i'} \to \{0,1\}\times\{0,1\}$ as 

\EQ{
f(m) := (\biggr[p \in \Ker\ a(m)\biggr], \biggr[q \in \Ker\ a(m)\biggr])
}
In the sequence $(f(i))_{i\geq i'}$, there cannot be any occurrence of $(0,0)$ because $f(j)=(0,0)$ implies $\Char_k(j)=0$.  There also cannot be any occurrences of $(1,1)$ because $f(j)=(1,1)$ implies $f(j+2)=(0,0)$. The patterns $(1,0), (0,1), (1,0)$ and $(0,1), (1,0), (0,1)$ also cannot occur in the sequence $(f(i))_{i\geq i'}$ because they violate property (ii) of Theorem 1. So, the only possibility is that $(f(i))_{j\geq i'}$ follows the pattern $(1,0),(1,0),(0,1),(0,1)$ modulo 4 (after some suitable offset.)

Recall that whenever a new bit $2^{a}$ appears in the sequence for the first time, it must occur in a number of the form $2^{a} + 2^{b}$. If a new bit were introduced in position $0$ or position $2$ of the pattern  $(1,0),(1,0),(0,1),(0,1)$, it would be impossible for the support of the term where the new bit is introduced to intersect with the support of the preceding term. Thus, new bits can only be introduced in positions $1$ or $3$. Here is a diagram illustrating what happens when a new bit $B$ is introduced in position $1$ (a blank indicates that bit is not specified, a 1 indicates that column's bit is present in the support, a 0 indicates it is absent:)
\begingroup
\setlength\arraycolsep{0pt}
\[
\begin{matrix*}[l]
  &\  &pq &\cdots     &B\\
  &0:\ &10 &\ &\ \\ 
  &1:\ &10000&\cdots0000&1\\
  &2:\ &01    &\         &1\\
  &3:\ &01    &\         & \\
\end{matrix*}
\]
\endgroup

Let $m$ be any positive integer not appearing in the first $i'$ terms of the sequence such that $\pop(m)\geq 2$ and $\Ker(m)\cap \Ker(k)=\emptyset$. Out of the infinitely many terms where a new bit $2^B$ is introduced, it is almost always the case that $B \notin \Ker(m)$. This implies that there are an infinite number of $0$s in the truth table of $\Char_m$. By Lemma 1, the number of $1$s occurring in the truth table of $\Char_m$ is also infinite. Thus, there are an infinite number of occurrences of $01$ in the truth table of $\Char_m$. By Lemma 2, this implies that $m$ appears in the sequence at some point after the $i'$th term. Since $\Ker(m)$ is disjoint from $\Ker(k)$, this is a contradiction of the assumption that $\Char_k(i)=1$ for all $i > i'$. 
\end{proof}

\begin{theorem6} 
 Every positive integer $n$ with $\pop(n)\geq 2$ appears in the sequence.
\end{theorem6}
\begin{proof}
Let $k$ be a positive integer with $\pop(k)=2$. An implication of Lemma 1 is that there is no positive integer $i'$ such that $i > i'$ implies $\Char_k(i) = 0$. An implication of Lemma 3 is that there is no $i'$ such that $i > i'$ implies $\Char_k(i)=1$. This proves that there there are infinite occurrences of $01$ in the truth table of $\Char_k$. By Theorem 4, the term $k$ must appear in the sequence.

By the previous paragraph, for any positive integer $n$ with $\pop(n)\geq 2$ there are infinitely many positive integers $x$ with $\pop(x)=2$ and $\Ker(x)\ \cap\ \Ker(n)\neq \emptyset$ appearing in the sequence. There are also infinitely many positive integers $y$ with $\pop(y)=2$  and $\Ker(y)\cap \Ker(n)= \emptyset$ appearing in the sequence. This implies that there are infinitely many occurrences of $01$ in the truth table of $\Char(n)$. By Theorem 4, $n$ must appear in the sequence.
\end{proof}

\medskip
\bibliographystyle{jis}
\bibliography{main_enots}

\bigskip
\hrule
\bigskip
2020 \emph{Mathematics Subject Classification}:~Primary 11B83. 
Secondary 11B75.
\medskip

\noindent 
\emph{Keywords}:~Enots Wolley,
Yellowstone Permutation,
permutation of the nonpowers of 2.
\bigskip
\hrule
\bigskip
\noindent 
(Concerned with sequences \seqnum{A098550}, \seqnum{A336957}, \seqnum{A338833}.)

\end{document}